\documentclass[preprint,12pt]{elsarticle}

\usepackage{epsfig}
\usepackage{amsmath,amssymb}
\usepackage{amsthm}
\usepackage{algorithm,algorithmic}
\usepackage{color,latexsym,bm}

\usepackage{hyperref}
\begin{document}

\begin{frontmatter}

\title{Preservation of the invariants of Lotka-Volterra equations by iterated deferred correction methods}

\author[sin]{Murat Uzunca\corref{cor1}}
\ead{muzunca@sinop.edu.tr}

\cortext[cor1]{Corresponding author}

\address[sin]{Department of Mathematics, Sinop University, 57000, Sinop, Turkey}

\begin{abstract}

In this paper we apply Kahan's nonstandard discretization to  three dimensional Lotka-Volterra equations in bi-Hamiltonian form. The periodicity of the solutions and all polynomial and non-polynomial invariants are well preserved in long-term integration. Applying classical deferred correction method, we show that the invariants are preserved with increasing accuracy as a results of more accurate numerical solutions. Substantial speedups over the Kahan's method are achieved at each run with deferred correction method.

\end{abstract}


\end{frontmatter}

%
%%%%%%%%%%%%%%%%%%%%%%%%%%%%%%%%%%%%%%%%%%%%%%%%%%%%%%%%%%%%%%%%%%%%%%%%%%%%%%%%%%%%%
%%%%%%%%%%%%%%%%%%%%%%%%%%%%%%%%%%%%%%%%%%%%%%%%%%%%%%%%%%%%%%%%%%%%%%%%%%%%%%%%%%%%%
\section{Introduction}

In the last two decades, many structure preserving geometric integrators are developed to preserve
symplectic structure, energy and other invariants,  phase space volume, reversing
symmetries, dissipation approximately or exactly (up to the round-off errors) \cite{Hairer10,Mclachlan98} of dynamical systems. These are symplectic and variational integrators for Hamiltonian systems  \cite{Hairer10,Marsden01}, integral preserving methods \cite{Mclachlan98}  and discrete  gradient methods \cite{Quispel96}. For special classes of
ordinary differential equations (ODEs), there exist non-standard discretization methods \cite{Kahan93,Mickens03} which preserve the conserved quantities and other features approximately or exactly. Among them Kahan’s method, also known as Hirota-Kimura method, applied to ODEs with quadratic vector fields,
preserves the integrals or conserved quantities of many Hamiltonian and integrable
systems \cite{Celledoni09,Celledoni14,Kimura00,Petrera07,Petrera11,Petrera17}. It was introduced by W. Kahan as "unconventional" discretization method  \cite{Kahan93} for quadratic vector fields
and applied to a scalar Riccati equation and a two-dimensional Lotka-Volterra system \cite{Kahan97}.

The Lotka-Volterra systems (LVSs) are first order ODEs, initially designed as an ecological predator-prey model. They occur in epidemiology,  in laser physics  \cite{Lamb64}, in evolutionary  game theory \cite{Hofbauer13} and as spatial discretizations of the Korteweg de Vries equation \cite{Kahan97,Karasozen06}. Most of the two and three dimensional LVSs  have periodic solutions and posses polynomial and non-polynomial integrals. They can be written as Poisson systems in bi-Hamiltonian form \cite{Magri78} and Nambu systems \cite{Nambu73}. Many numerical methods are applied to LVSs which preserve the integrals, periodic solution, attractors and son on \cite{Mickens03,Roeger05,Roeger06,Ionescu15,Cohen11,SanzSerna94}.

For Hamiltonian systems, higher order accuracy for integrals can be  achieved by composing
symplectic integrators in time  \cite{Blanes08,McLachlan02,Suzuki90,Yoshida90}. Starting with a basic method, arbitrary orders of accuracy can be obtained by applying the composition  to a  lower order symplectic method recursively.    Another class of numerical methods designed
for the construction of high-order approximations to the solution of differential
equations are the deferred correction methods.
A numerical solution of an initial-value problem (IVP)
for ODEs is computed by a low order method and then subsequently
refined by solving the IVP constructed by the error between the numerical and continuous solutions. Under suitable assumptions, this process can
be repeated to produce solutions with an arbitrarily high order of accuracy.
Deferred correction methods have been extensively applied to IVPs
such as, classical deferred correction (CDC) methods \cite{Dutt00,Hansen11},  spectral deferred correction methods \cite{Dutt00,Minion03} and integral deferred correction methods  \cite{Christlieb10}.

The LVSs are ODEs with quadratic vector fields. In this paper, two three dimensional (3D) LVSs in bi-Hamiltonian form are solved by  the CDC method based on Kahan's method. We show that the periodicity of the solutions and integrals are preserved in long term integration. At each correction step, more accurate solutions are obtained and the integrals are preserved more accurately. Iterated deferred correction methods are more efficient than the composition methods, because at the correction step the same grid is used. Therefore substantial speedups can be obtained by the CDC methods over the basic method, i.e. Kahan's method.  To the best of our knowledge, the deferred correction methods are used first time to preserve the conserved quantities of dynamical systems with higher accuracy.

The paper is organized as follows. In Section~\ref{sec:lvs},  we present two 3D LVSs in bi-Hamiltonian form. In Section~\ref{sec:kahan}, we give a short description of Kahan's method applied to ODEs with a quadratic vector field.  Algorithm for the CDC methods is discussed and given in Section~\ref{sec:cdc}. Numerical results for Kahan's method with CDC methods and composition methods are compared in Section~\ref{sec:numeric}. The paper ends with some conclusions in Section~\ref{sec:conclusion}.

%%%%%%%%%%%%%%%%%%%%%%%%%%%%%%%%%%%%%%%%%%%%%%%%%%%%%%%%%%%%%%%%%%%%%%%%%%%%%%%%%%%%%
%%%%%%%%%%%%%%%%%%%%%%%%%%%%%%%%%%%%%%%%%%%%%%%%%%%%%%%%%%%%%%%%%%%%%%%%%%%%%%%%%%%%%
\section{Lotka-Volterra systems} 
\label{sec:lvs}

The LVSs are systems of first order ODEs in the following form \cite{Maier13,Schmming03}
\begin{equation}\label{lvs}
\dot{u}_i  = u_i(r_i + \sum_j \alpha_{ij} u_j ),\quad i=1,\ldots,m,
\end{equation}
where ${\mathbf u} := (u_1,\ldots u_m)^T$ is the $m$-dimensional state vector and  $\dot{u}_i =du_i/dt$ denotes the derivative with respect to time.
In ecology, $u_i$ describe the densities of each species and $r_i$ are the intrinsic growth or decay rates. The interaction between the species is specified by the coefficient matrix  $A = (\alpha_{ij}), \; i,j=1,\ldots, m$.  All variables in \eqref{lvs} are real and the densities $u_i$ are positive.
There are no closed solutions of LVSs when $m\ge 2$, they have to be integrated numerically.

%%%%%%%%%%%%%%%%%%%%%%%%%%%%%%%%%%%%%%%%%%%%%%%%%%%%%%%%%%%%%%%%%%%%%%%%%%%%%%%
\subsection{Bi-Hamiltonian 3D Lotka-Volterra systems}

Many 3D LVSs can be written in the following bi-Hamiltonian form

\begin{equation}\label{biham}
\dot{{\mathbf u}} = J_1\nabla H_2 = J_2\nabla H_1,
\end{equation}
where $J_1$ and $J_2$ are the skew-symmetric Poisson matrices satisfying the Jacobi identity.
There exists two
independent integrals $H_1$ and $H_2$, associated with $J_1$ and $J_2$ such
that $H_2$ is the Casimir for one Poisson structure while $H_1$ is the Casimir for the
other \cite{Nutku90}.
Bi-Hamiltonian systems are completely integrable \cite{Magri78}.
3D LVSs can also be written as Nambu systems \cite{Nambu73}, as generalization of Hamiltonian systems with multiple Hamiltonians. Nambu form of
\eqref{biham} is given as
\begin{equation*}
\dot{\mathbf u} = \nabla H_1 \times  \nabla H_2.
\end{equation*}
Vector fields of Nambu systems are divergence free and the flow is volume preserving.

A well known 3D LVS possessing bi-Hamiltonian structure is \cite{Grammaticos90,Ionescu15,Plank96}
\begin{equation}\label{lv1}
\begin{aligned}
\dot{u}_1 & = u_1(cu_2+u_3 + \lambda), \\
\dot{u}_2 & = u_2(u_1+au_3 + \mu), \\
\dot{u}_3 & = u_3(bu_1+u_2 + \nu),
\end{aligned}
\end{equation}
where $\lambda,\mu,\nu > 0$, and with $abc=-1$ and $\nu =\mu b -\lambda ab$.
The skew-symmetric Poisson matrices and Hamiltonians then are given by
$$
J_1 = \left ( \begin{array}{ccc}
 0 & cu_1u_2 & bc u_1u_3 \\
 -cu_1u_2 & 0 & -u_2u_3 \\
 -bcu_1u_3 & u_2u_3 & 0
 \end{array} \right ),
$$
$$
  J_2 = \left ( \begin{array}{ccc}
 0 & cu_1u_2(au_3 +\mu) &  c u_1u_3(u_2 +\nu) \\
 -cu_1u_2(au_3 +\mu) & 0 & u_1u_2u_3 \\
  -c u_1u_3(u_2 +\nu) & -u_1u_2u_3 & 0
 \end{array} \right ),
$$
$$H_1 = ab\ln u_1 -b \ln u_2 + \ln u_3,\;  H_2=abu_1 +u_2 -a u_3 +\nu\ln u_2 -\mu \ln u_3.
$$
$H_1$ and $H_2$ are Casimirs of $J_1$ and $J_2$, respectively, i.e. $J_1\nabla H_1=0$ and $J_2\nabla H_2 = 0$.

Another example of 3D LVS  is  the reversible 3D LVS  with the circulant coefficient matrix $A$ \cite{Yang12}
\begin{equation}\label{lvs2}
\begin{aligned}
\dot{u}_1 & = u_1(u_2-u_3), \\
\dot{u}_2 & = u_2(u_3 -u_1), \\
\dot{u}_3 & = u_3(u_1 -u_2).
\end{aligned}
\end{equation}
It has a game-theoretical interpretation \cite{Hofbauer13} and possesses bi-Hamiltonian form
with the Poisson matrices
$$
J_1 = \left ( \begin{array}{ccc}
0 & -1 & 1 \\
 1 & 0 & -1 \\
 -1 & 1 & 0
 \end{array} \right ),\quad
  J_2 = \left ( \begin{array}{ccc}
 0 & u_1u_2 &  -u_1u_3 \\
 -u_1u_2 & 0 & u_2u_3 \\
 u_1u_3 & -u_2u_3 & 0
 \end{array} \right ),
$$
and with the linear Hamiltonian $H_1 = u_1 + u_2 + u_3$ and the cubic Hamiltonian $H_2=u_1u_2u_3$.
It is reversible with respect to $\rho = {\mathrm diag}(-1,-1,-1)$, $\rho f(\mathbf{x}) = - f(\rho{\mathbf x})$. It can also be written as Nambu system. The flow generated by \eqref{lvs2} is source free, i.e. the volume is preserved. The linear integral $H_1$ represents the volume. The $n$-dimensional extension of \eqref{lvs2} as integrable discretization of the Korteweg de Vries equation was integrated with a Poisson structure preserving integrator in \cite{Karasozen06}.
Necessary and sufficient conditions for conservation laws of  $n$-dimensional LVSs \eqref{lvs} including the two Poisson systems are derived in
\cite{Schmming03}.

%%%%%%%%%%%%%%%%%%%%%%%%%%%%%%%%%%%%%%%%%%%%%%%%%%%%%%%%%%%%%%%%%%%%%%%%%%%%%%%
%%%%%%%%%%%%%%%%%%%%%%%%%%%%%%%%%%%%%%%%%%%%%%%%%%%%%%%%%%%%%%%%%%%%%%%%%%%%%%%
\section{Kahan's method }
\label{sec:kahan}

The LVS \eqref{lvs} is an autonomous ODE system in the following form
\begin{equation}\label{quad}
\dot{\mathbf u} = f({\mathbf u}) := Q({\mathbf u})+B{\mathbf u},
\end{equation}
with the quadratic vector field  $(Q({\mathbf u}))_i = u_i\sum_j \alpha_{ij} u_j$ and the diagonal matrix $B = diag(r_1,\ldots , r_m)$. In the system \eqref{quad}, the unknown solution vector is ${\mathbf u}=(u_1,\ldots , u_m)^T$, and it is prescribed the vector of initial conditions ${\mathbf u}(t_0) = (u_1(t_0),\ldots , u_m(t_0))^T$.

For the ODE system \eqref{quad}, Kahan introduced in 1993 the "unconventional"  discretization as \cite{Kahan93}
\begin{equation*}
\frac{{\mathbf u}_{n+1} - {\mathbf u}_n}{\Delta t} = Q({\mathbf u}_n,{\mathbf u}_{n+1}) + \frac{1}{2}B({\mathbf u}_n + {\mathbf u}_{n+1}) ,
\end{equation*}
where $\Delta t $ is the step size of the integration, ${\mathbf u}_{n+1}$ and   ${\mathbf u}_n$ are the approximations at the time instances $t_{n+1}$ and $t_n$, respectively. The symmetric bilinear form $Q(\cdot ,\cdot )$ is obtained by the polarization of the quadratic vector field $Q(\cdot)$ \cite{Celledoni15}
$$
Q({\mathbf u}_n,{\mathbf u}_{n+1}) = \frac{1}{2}\left( Q({\mathbf u}_n+{\mathbf u}_{n+1}) - Q({\mathbf u}_n) -  Q({\mathbf u}_{n+1}) \right).
$$

The Kahan's method is second order and time-reversal \cite{Kahan97}:
\begin{align*}
\quad \frac{{\mathbf u}_{n+1} - {\mathbf u}_n}{\Delta t} & = \left( I -\frac{\Delta t}{2} f'({\mathbf u}_n)\right)^{-1} f({\mathbf u}_n), \\
\quad \frac{{\mathbf u}_{n+1} - {\mathbf u}_n}{\Delta t} & =  \left ( I +\frac{\Delta t}{2} f'({\mathbf u}_{n+1})\right )^{-1} f({\mathbf u}_{n+1}),
\end{align*}
where $I\in \mathbb{R}^{m\times m}$ is the identity matrix and $f'$ denotes the Jacobian of $f$.
Moreover, Kahan's method is linearly implicit and it coincides with a certain
Rosenbrock method on quadratic vector fields, i.e. ${\mathbf u}_{n+1}$ can be computed by solving a single linear system
$$
\left ( I -\frac{\Delta t}{2} f'({\mathbf u}_n)\right ) \tilde {\mathbf u} = \Delta t f({\mathbf u}_n),\qquad {\mathbf u}_{n+1} = {\mathbf u}_n + \tilde {\mathbf u}.
$$
Symplectic integrators like the implicit mid-point rule \cite{Hairer10}, energy preserving average vector field method \cite{Cohen11} and conservative methods \cite{Wan17} require at each time step more than one Newton iteration to solve nonlinear implicit equations to preserve the integrals accurately. Due to the linearly implicit nature,  Kahan's method is a very efficient structure preserving integrator for ODEs with quadratic vector fields.

Kahan's method is also a Runge-Kutta method, with negative weights, restricted to quadratic vector fields \cite{Celledoni13}:
\begin{equation*}
\frac{{\mathbf u}_{n+1} - {\mathbf u}_n}{\Delta t} = -\frac{1}{2}f({\mathbf u}_n) + 2f \left (\frac{{\mathbf u}_{n+1} + {\mathbf u}_n}{2}      \right )   -\frac{1}{2}f({\mathbf u}_{n+1}).
\end{equation*}

Kahan’s method was independently rediscovered by Hirota and Kimura \cite{Hirota00,Kimura00}, which preserves the integrability for  a large number of integrable quadratic
vector fields like Euler top, Lagrange top \cite{Kimura00,Petrera07,Petrera17}  Suslov and Ishii systems, Nambu systems, Riccati systems,
and the first Painlev\'{e} equation \cite{Celledoni13,Celledoni15}.
Kahan’s method is generalized in \cite{Celledoni15} to cubic and  higher
degree polynomial vector fields.

%%%%%%%%%%%%%%%%%%%%%%%%%%%%%%%%%%%%%%%%%%%%%%%%%%%%%%%%%%%%%%%%%%%%%%%%%%%%%%%
%%%%%%%%%%%%%%%%%%%%%%%%%%%%%%%%%%%%%%%%%%%%%%%%%%%%%%%%%%%%%%%%%%%%%%%%%%%%%%%
\section{Iterative deferred correction method}
\label{sec:cdc}

In this section, we apply the CDC method \cite{Dutt00,Hansen11} to the ODE system \eqref{quad} with quadratic vector field, related to the LVSs \eqref{lvs}.
For a given time interval $[0,T]$, we subdivide it into $J$ equidistant intervals $[t_j,t_{j+1}]$, $j=0,1,\ldots , J-1$:
$$
0=t_0 < t_1 < \cdots < t_j < \cdots < t_J=T, \quad t_{j+1}-t_j = \Delta t,
$$
on which we define the approximate solutions by ${\mathbf u}_{j}\approx {\mathbf u}(t_j)$, $j=1,\ldots , J$, and for $j=0$ we use the initial condition, ${\mathbf u}_{0}= {\mathbf u}(0)$. Further, each interval $[t_j,t_{j+1}]$ is subdivided into $n-1$ equidistant intervals forming $n$ nodes including the end points $t_j$ and $t_{j+1}$ as
$$
t_j=t_{j,1} < t_{j,2} < \cdots < t_{j,i} < \cdots < t_{j,n} =t_{j+1},
$$
and we define the approximate solutions on these nodes by ${\mathbf u}_{j,i}\approx {\mathbf u}(t_{j,i})$, $i=1,\ldots , n$, given that ${\mathbf u}_{j,1}= {\mathbf u}_j$.
The CDC method, on each subinterval  $[t_j,t_{j+1}]$, starts by solving the ODE system \eqref{quad} for the solutions at the nodes $\{t_{j,i}\}_{i=1}^n$, with a method of order $p_0$. Then, the approximate solutions of the ODE system \eqref{quad} on the interval $[t_j,t_{j+1}]$ are defined by ${\mathbf U}_j^{[0]} := ({\mathbf u}^{[0]}_{j,1},\ldots,
{\mathbf u}^{[0]}_{j,n} )$, and satisfy that
$$
{\mathbf u}_{j,i} ^{[0]} = {\mathbf u}(t_{j,i}) + {\mathcal O} ((\Delta t)^{p_0}), \quad i=1,\ldots ,n.
$$
We apply here the second order Kahan's method described in Section~\ref{sec:kahan} as the basic method, so $p_0=2$.
After, it follows the correction procedure. At the $s$-th correction step, the CDC method computes an improved (corrected) solution ${\mathbf U}_j^{[s]} := ({\mathbf u}^{[s]}_{j,1},\ldots, {\mathbf u}^{[s]}_{j,n} )$ of the following error system
\begin{equation}\label{errorsystem}
\begin{aligned}
\frac{d}{dt} \mathbf{e}_{j}^{[s-1]}(t,\mathbf{U}_j^{[s-1]}(t)) &= f\left(\mathbf{e}_{j}^{[s-1]}(t,\mathbf{U}_j^{[s-1]}(t)) + \mathbf{U}_j^{[s-1]}(t) \right)  - \frac{d}{dt} \mathbf{U}_j^{[s-1]}(t),\\
\mathbf{e}_{j}^{[s-1]}(t_j,\mathbf{U}_j^{[s-1]}(t_j)) &= 0,
\end{aligned}
\end{equation}
by a method of order $p_s$. In \eqref{errorsystem}, $\mathbf{e}_{j}^{[s]}(t,\mathbf{U}_j^{[s]}(t))$ denotes the error function on the $s$-th iteration step given by
\begin{equation}\label{errors}
{\mathbf e}_j^{[s]}(t,{\mathbf U}_j^{[s]}(t)) = {\mathbf u}(t) -  {\mathbf U}_j^{[s]}(t), \quad t \in [t_j,t_{j+1}].
\end{equation}
The differences between the deferred correction methods are based on the formation of an error system; on the continuous level they are equivalent. In the CDC method, the error system \eqref{errorsystem} is used, which is obtained by the differentiation of the error equation \eqref{errors} with respect to the time variable $t$. The function ${\mathbf U}_j^{[s]}(t)$ stands for the continuous approximation of the discrete solutions ${\mathbf U}_j^{[s]} = ({\mathbf u}^{[s]}_{j,1},\ldots, {\mathbf u}^{[s]}_{j,n} )$.
Here, we construct the continuous approximation ${\mathbf U}_j^{[s]}(t)$ based on the Lagrange interpolation as
$$
\mathbf{U}_j^{[s]}(t) = \sum_{k=1}^n l_k(t)\cdot \mathbf{u}_{j,k}^{[s]}, \quad
l_k(t) = \prod_{i\neq k}\frac{t - t_{j,i}}{t_{j,k}-t_{j,i}},
$$
where $l_k(t)$ are the Lagrange basis functions.

The error system  \eqref{errorsystem} is non-autonomous due the occurrence of the time dependent terms $\mathbf{U}_j^{[s-1]}(t)$ and their derivatives.   Because Kahan's method is designed for autonomous systems, in the correction steps we use the second order  mid-point rule.

After, defining the vector of error approximations ${\mathbf E}_j^{[s]} := ({\mathbf e}^{[s]}_{j,1},\ldots, {\mathbf e}^{[s]}_{j,n} )$ where ${\mathbf e}^{[s]}_{j,i}$ are the discrete solutions of the error system \eqref{errorsystem} on the nodes $\{t_{j,i}\}$, we obtain the corrected numerical approximations through the update formula
$$
\mathbf{U}_j^{[s]} = \mathbf{U}_j^{[s-1]} + \mathbf{E}_j^{[s-1]}.
$$
An outline of the CDC method can be found in Algorithm~\ref{alg:cdc}.

\begin{algorithm}[t]
\caption{Classical deferred correction method \label{alg:cdc}}
\textbf{Input:} Correction number $S$, partition $\{[t_j,t_{j+1}]\}_{j=1}^J$ of the time interval $[0,T]$, initial solution $\mathbf{u}_0:= \mathbf{u}(0)$
 \\
\textbf{Output:} The approximate solutions $\{\mathbf{u}_1 , \ldots , \mathbf{u}_J\}$
 \\
\begin{algorithmic}[1]
\FOR{$j=0,1,\ldots , J-1$}
	\STATE Set $\mathbf{u}_{j,1}:=\mathbf{u}_j$
	\STATE Solve the ODE system \eqref{quad} for ${\mathbf U}_j^{[0]}= ({\mathbf u}_{j,1}^{[0]},\ldots, {\mathbf u}_{j,n}^{[0]} )$ on the nodes $\{t_{j,i}\}_{i=1}^n$, using Kahan's method\\
	\FOR{$s=1,2,\ldots , S$}
		\STATE Form the continuous solution $\mathbf{U}_j^{[s-1]}(t)$ using discrete set $\mathbf{U}_j^{[s-1]}$
		\STATE Solve the error system \eqref{errorsystem} for ${\mathbf E}_j^{[s-1]} = ({\mathbf e}_{j,1}^{[s-1]},\ldots, {\mathbf e}_{j,n}^{[s-1]} )$ on the nodes $\{t_{j,i}\}_{i=1}^n$, using mid-point rule
		\STATE Update the solution vector as $\mathbf{U}_j^{[s]} = \mathbf{U}_j^{[s-1]} + \mathbf{E}_j^{[s-1]}$
	\ENDFOR
	\STATE Set the solution $\mathbf{u}_{j+1}:= {\mathbf u}_{j,n}^{[S]}$
\ENDFOR
\end{algorithmic}
\end{algorithm}

Expected order of accuracy of the CDC methods for uniformly spaced nodes is given by $\min \{ P_S , n-1\}$, where $P_S = \sum_{s=0}^S p_s$, $S$ is the number of corrections and $n$ is the number of nodes used in each interval $[t_j,t_{j+1}]$ \cite{Dutt00,Hansen11}. Since we use Kahan's method and mid-point method, both of which are second order methods, we have $p_s=2$ for all $s=0,1,\ldots ,S$, and then the expected order of accuracy becomes $\min \{ 2S+2,n-1\}$. According to this fact, we set $n=2S+3$ in the simulations to obtain the expected order of accuracy as $n-1$.
When non-uniform nodes  like Gauss--Lobatto, Gauss--Legendre, and Chebyshev nodes are used,  for CDC methods
the accuracy improves with more corrections although the order
of accuracy stagnates at two  \cite{Hansen11}. When a low order Lagrange interpolation is used on small intervals $[t_j,t_{j+1}]$, the CDC method can produce accurate results on uniform nodes, as it will be shown in the next Section.

%%%%%%%%%%%%%%%%%%%%%%%%%%%%%%%%%%%%%%%%%%%%%%%%%%%%%%%%%%%%%%%%%%%%%%%%%%%%%%%
%%%%%%%%%%%%%%%%%%%%%%%%%%%%%%%%%%%%%%%%%%%%%%%%%%%%%%%%%%%%%%%%%%%%%%%%%%%%%%%
\section{Numerical results}
\label{sec:numeric}

In this section, we present numerical results for the Lotka-Volterra systems described in Section \ref{sec:lvs} solving by Kahan's method, and demonstrate the performance of the CDC method.
In all examples, we give the results of the runs using Kahan's method with a small time step-size without CDC method, and the ones with the CDC method using Kahan's method for the ODE system \eqref{quad} and the mid-point rule for the error system \eqref{errorsystem}, with a larger time step-size.
Hamiltonian errors  $H(0)-H(t)$ are plotted over $t$. We set $S=1$ and accordingly $n=2S+3=5$ in the CDC procedure.

The $L^2$-error for a Hamiltonian $H(t)$, and the $L^2$-error between the exact solution $\mathbf{u}_{exact}(t)$ and the numerical solution $\mathbf{u}$ are measured using the following norms
$$
L^2(H) = \left( \Delta t\sum_{i=1}^J [H(t_i)-H(0)]^2 \right)^{1/2}\; , \quad L^2(\mathbf{u}) = \left( \Delta t\sum_{i=1}^J [\mathbf{u}_i-\mathbf{u}_{exact}(t_i)]^2 \right)^{1/2},
$$
where the exact solution $\mathbf{u}_{exact}(t)$ is obtained by MatLab's {\em ode45} solver in which we set the relative and the absolute tolerances as $10^{-13}$. The  order of accuracy is calculated as
$$
\text{order} = \frac{1}{\log 2}\log \left(  \frac{\text{Err}_{\Delta t_i}}{\text{Err}_{\Delta t_{i+1}}} \right),
$$
where $\text{Err}_{\Delta t_i}$ and $\text{Err}_{\Delta t_{i+1}}$ stand for the $L^2$-error of an Hamiltonian or the solution, obtained by the  consecutive step sizes $\Delta t_i$ and $\Delta t_{i+1}=\Delta t_i/2$, respectively.

%%%%%%%%%%%%%%%%%%%%%%%%%%%%%%%%%%%%%%%%%%%%%%%%%%%%%%%%%%%%%%%%%%%%%%%%%%%%%%%%%%%%%%%%
\subsection{Bi-Hamiltonian LVS}
\label{Ex1}

We consider the 3D LVS \eqref{lv1} on the interval $[0,100]$, with the parameter values
$
(a,b,c,\lambda , \mu , \nu)=(-1,-1,-1,0,1,-1)
$ \cite{Ionescu15}.
The initial condition is taken as $(u_1(0),u_2(0),u_3(0) )^T = ( 1 , 1.9 , 0.5 )^T$.

\begin{figure}[htb!]
\centering
\includegraphics[scale=0.4]{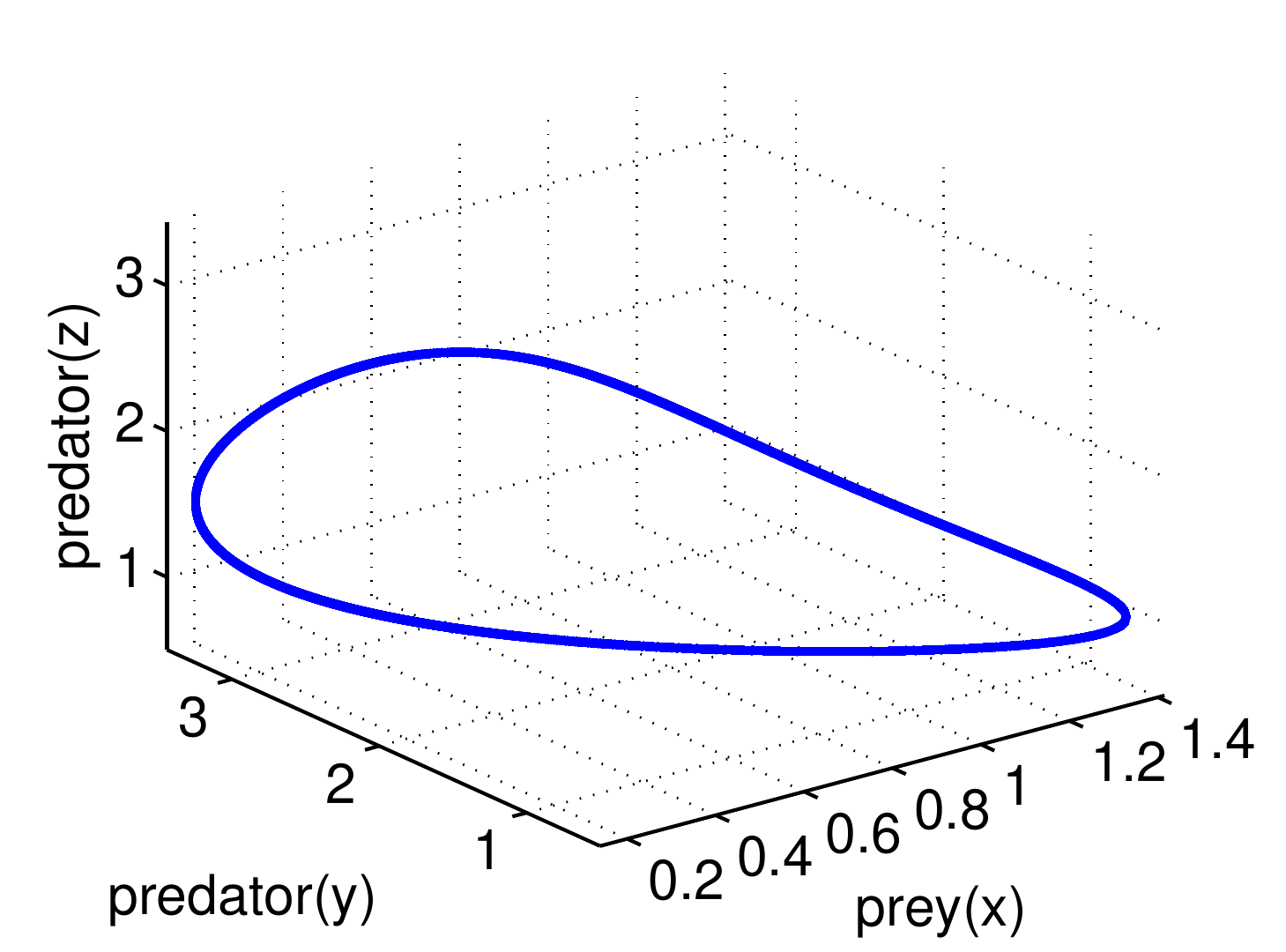}

\includegraphics[scale=0.4]{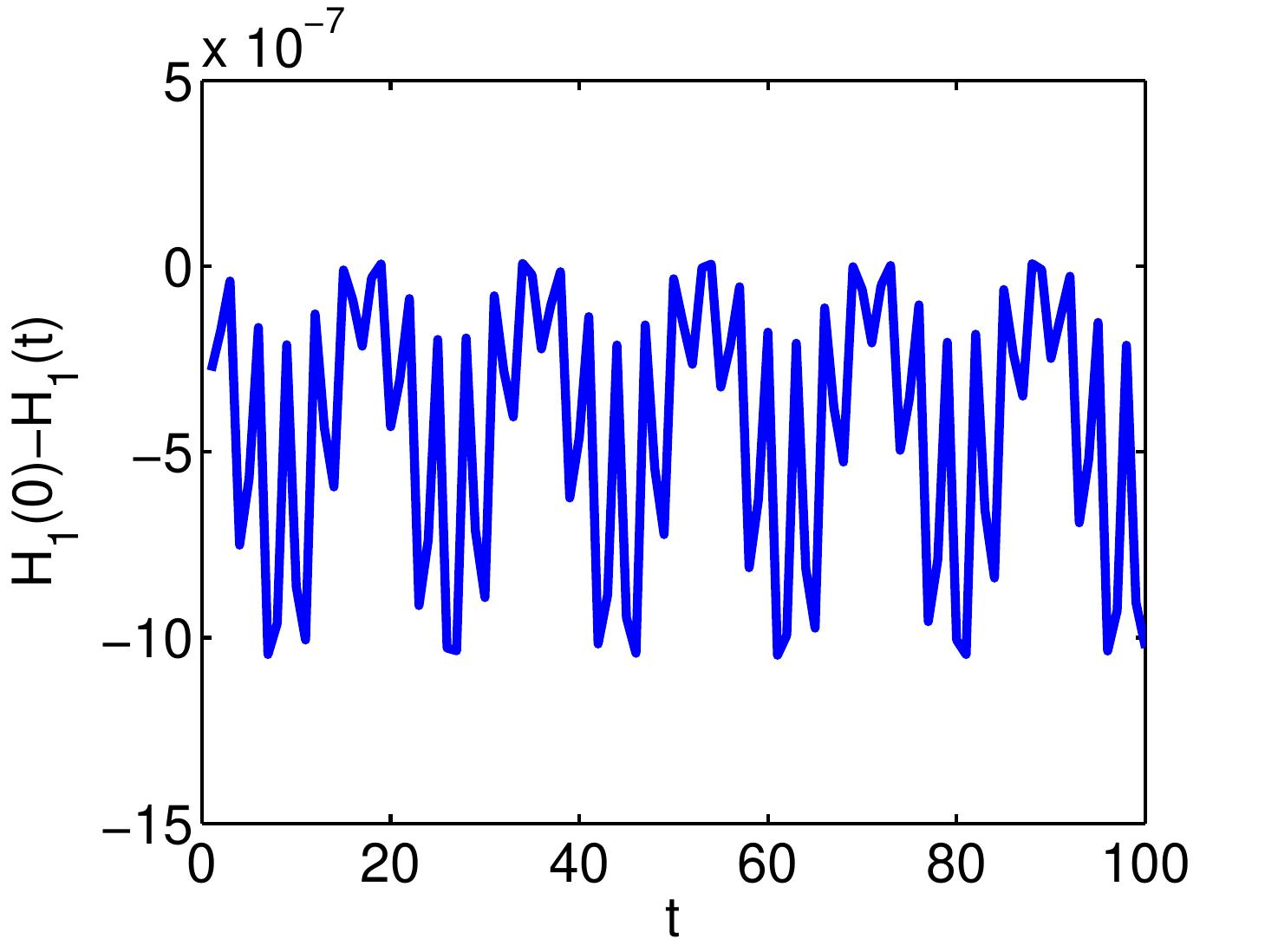}
\includegraphics[scale=0.4]{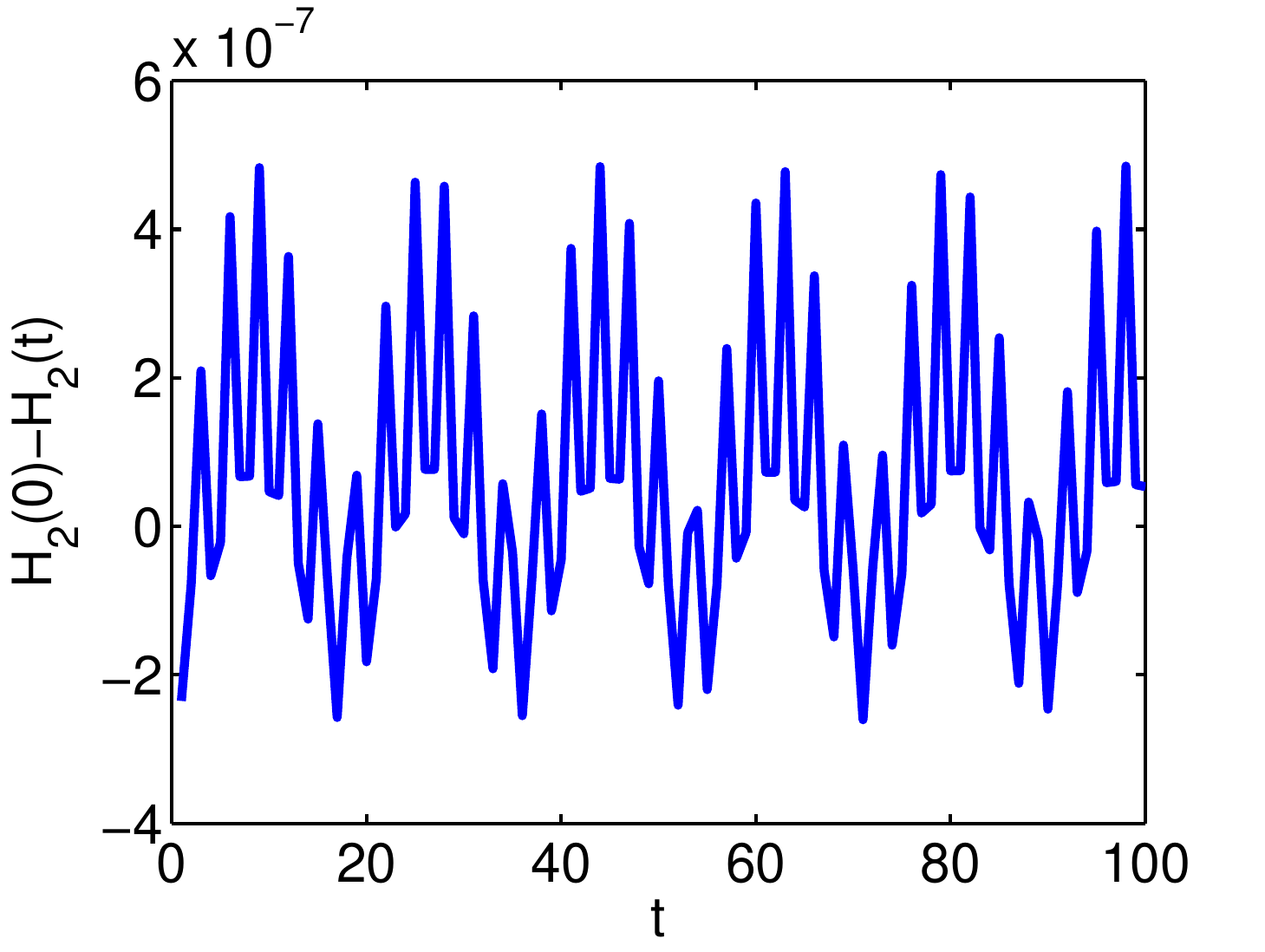}
\caption{Example~\ref{Ex1} without CDC: (Top) Periodic solutions; (Bottom) errors of Hamiltonians $H_1$ (left) and $H_2$ (right) from the initials: $\Delta t=0.001$}
\label{lv_fig_kahan}
\end{figure}

We show that Kahan's method preserves the periodicity of the solutions and the Hamiltonians in Figure~\ref{lv_fig_kahan}.
It was proved in \cite{Roeger05} that Kahan's method preserves the periodicity of LVSs \eqref{lvs}. The average vector field method, which preserves the Poisson structure, was also applied to LVSs \eqref{lvs} in \cite{Cohen11}. It was shown there that the first Hamiltonian $H_1$ of \eqref{lv1} is preserved, but the Casimir $H_2$ shows a drift in long term integration.

\begin{figure}[htb!]
\centering
\includegraphics[scale=0.4]{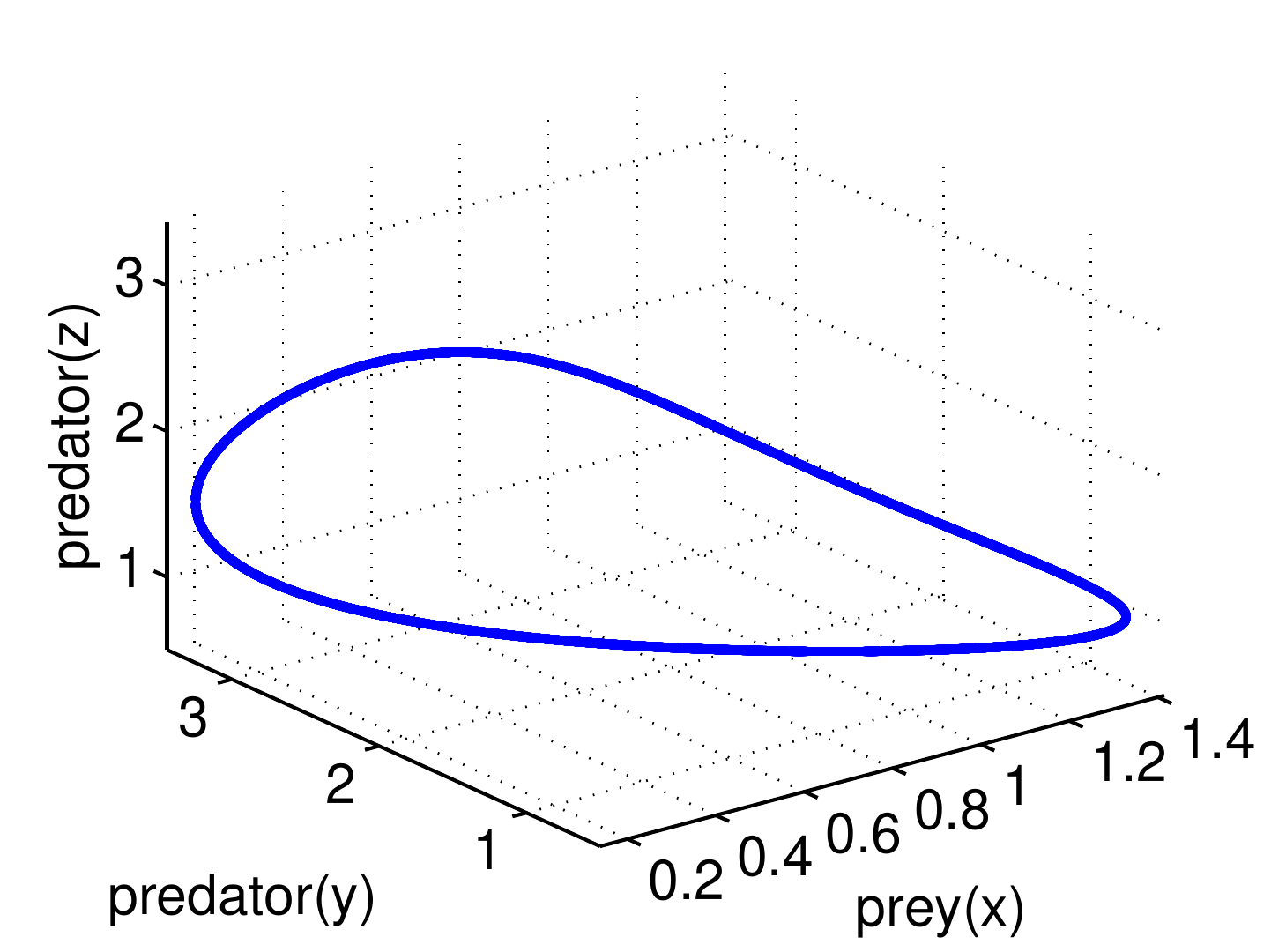}

\includegraphics[scale=0.4]{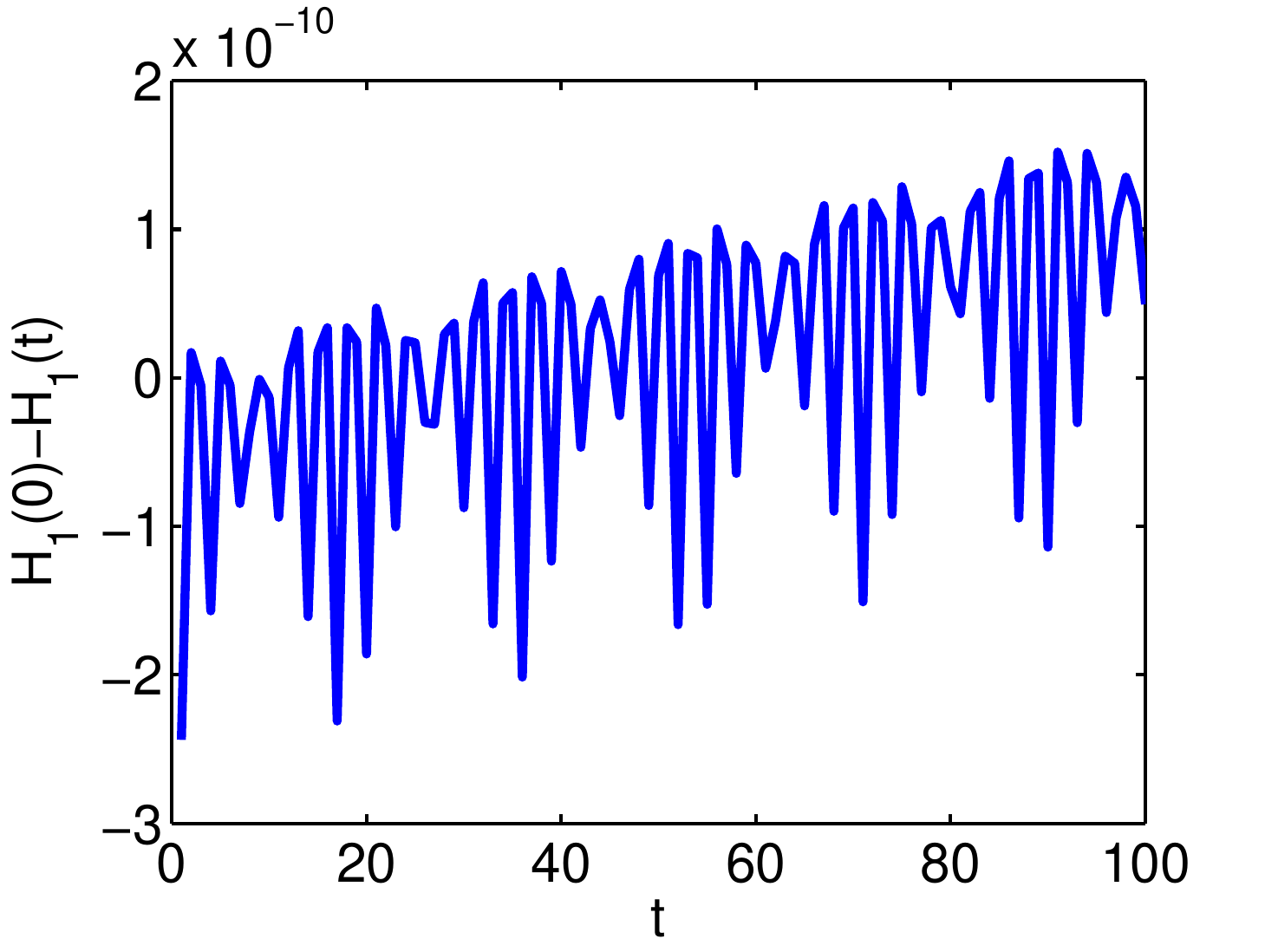}
\includegraphics[scale=0.4]{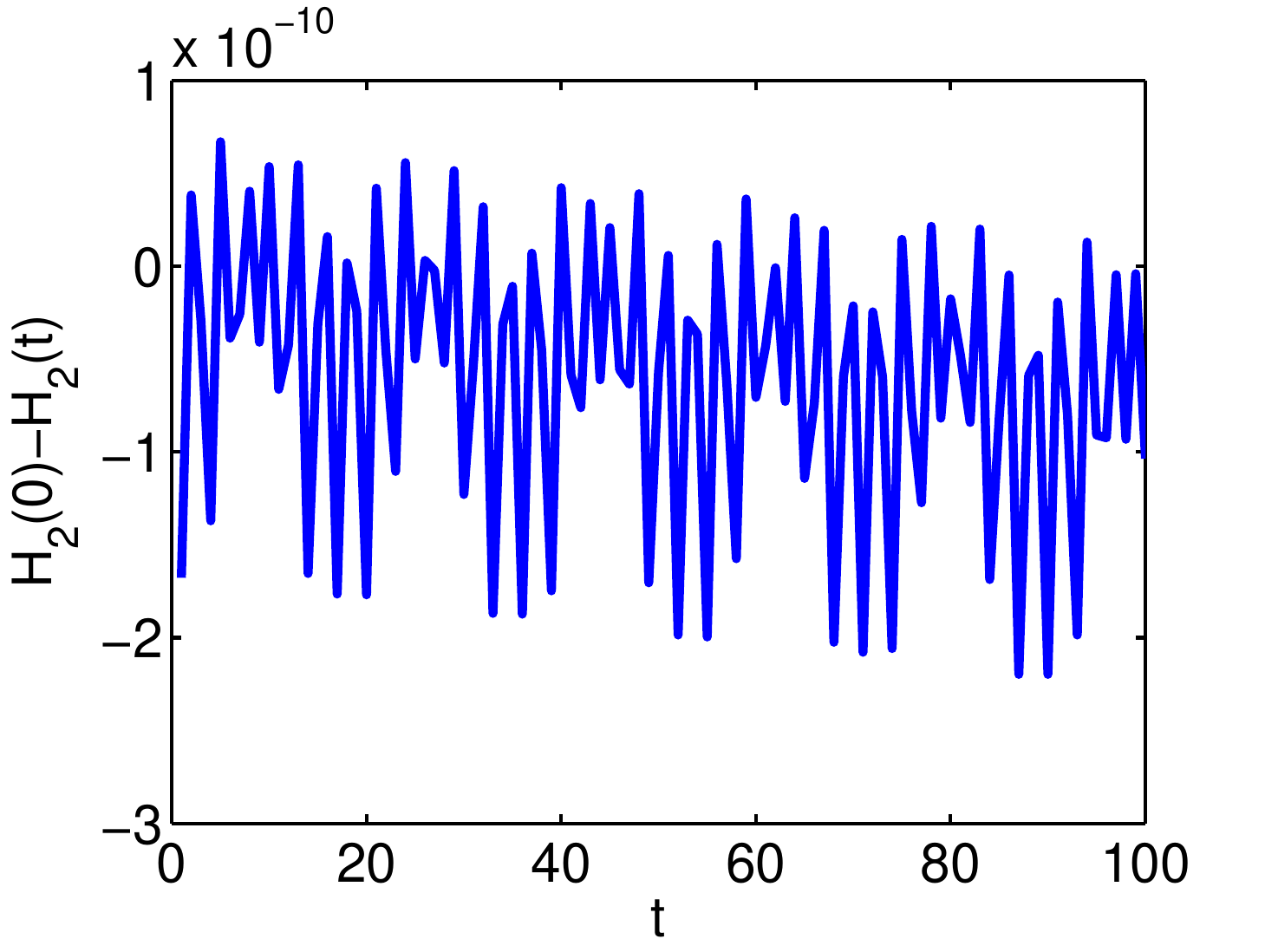}
\caption{Example~\ref{Ex1} with CDC: (Top) Periodic solutions; (Bottom) errors of Hamiltonians $H_1$ (left) and $H_2$ (right) from the initials: $\Delta t=0.01$, $S=1$, $n=5$}
\label{lv_fig_kahan_mid}
\end{figure}

When the CDC method is applied with $S=1$, $n=2S+3=5$ and with the  use of larger time step-size $\Delta t = 0.01$, the periodicity of the solutions and  Hamiltonians $H_1$ and $H_2$ are preserved  in  Figure~\ref{lv_fig_kahan_mid}.  Compared with the numerical results in Figure~\ref{lv_fig_kahan}, it turns out that the use of CDC method is more efficient in terms of the preservation of the Hamiltonians. In Figure~\ref{lv_fig_kahan_mid}, a slow drift in the preservation of the Hamiltonians is observed. When composition methods are applied to Kahan’s method, a comparatively more rapid Hamiltonian drift is observed \cite{Celledoni13}.

\begin{figure}[htb!]
\centering
\includegraphics[scale=0.26]{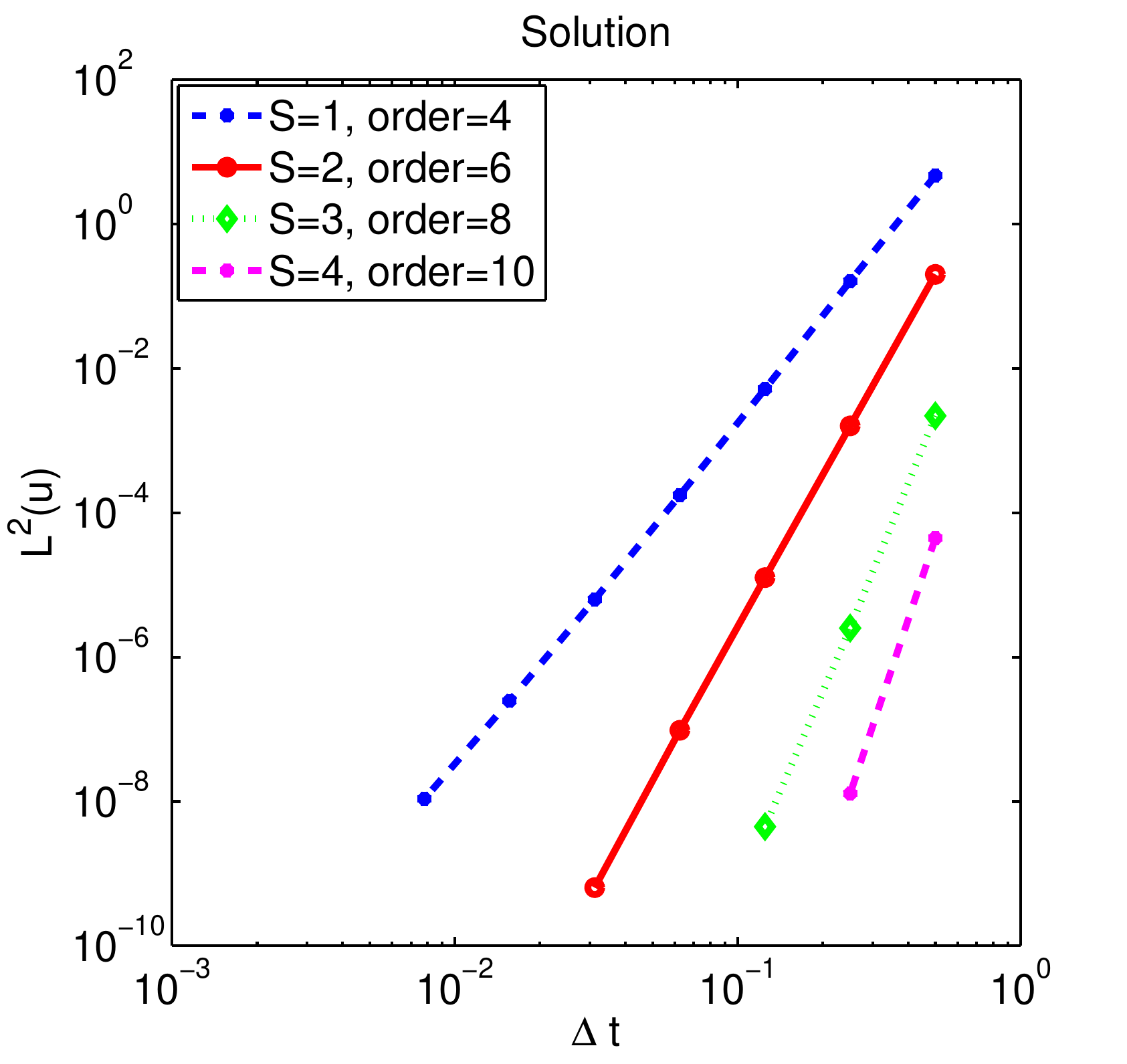}
\includegraphics[scale=0.26]{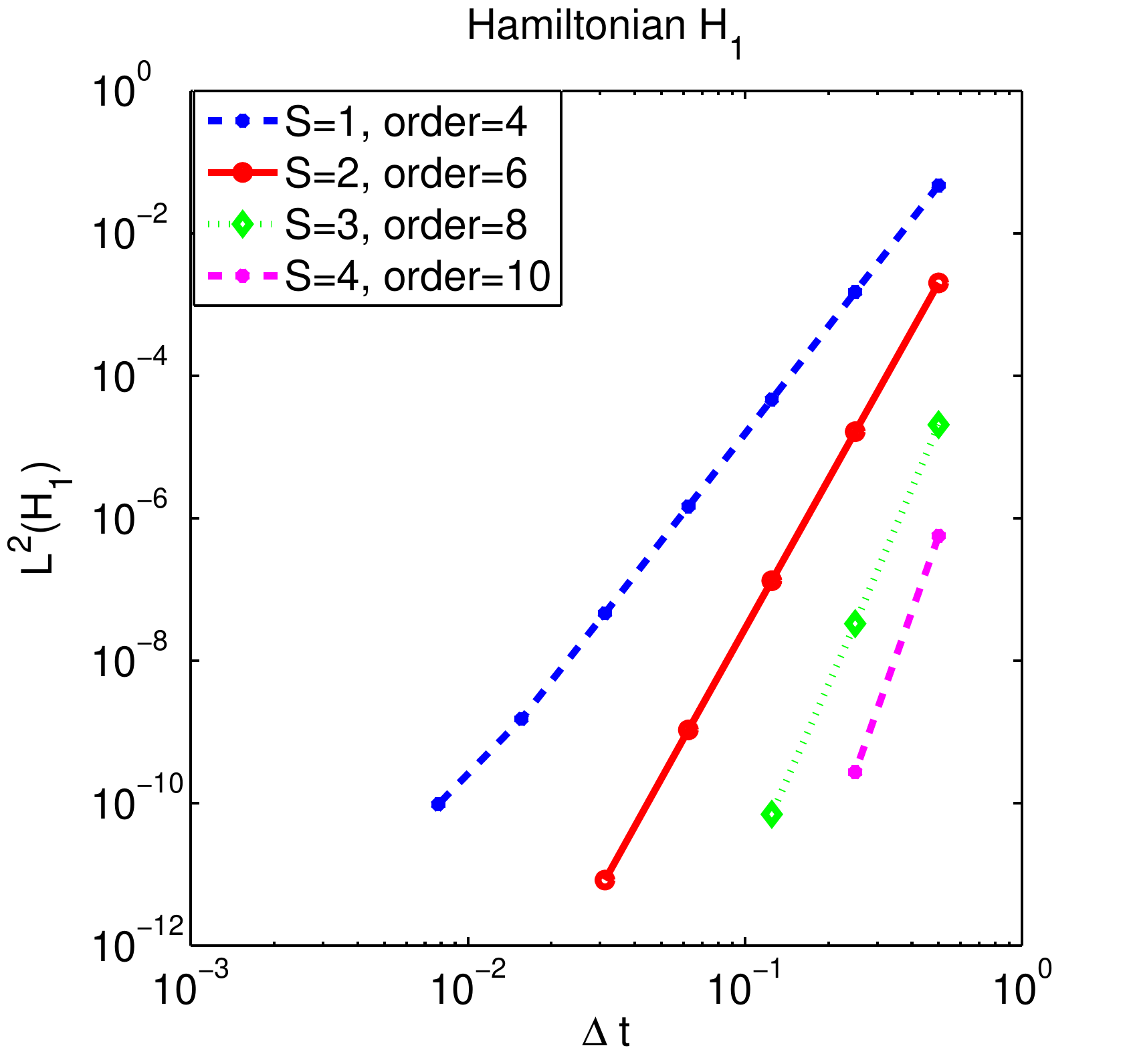}
\includegraphics[scale=0.26]{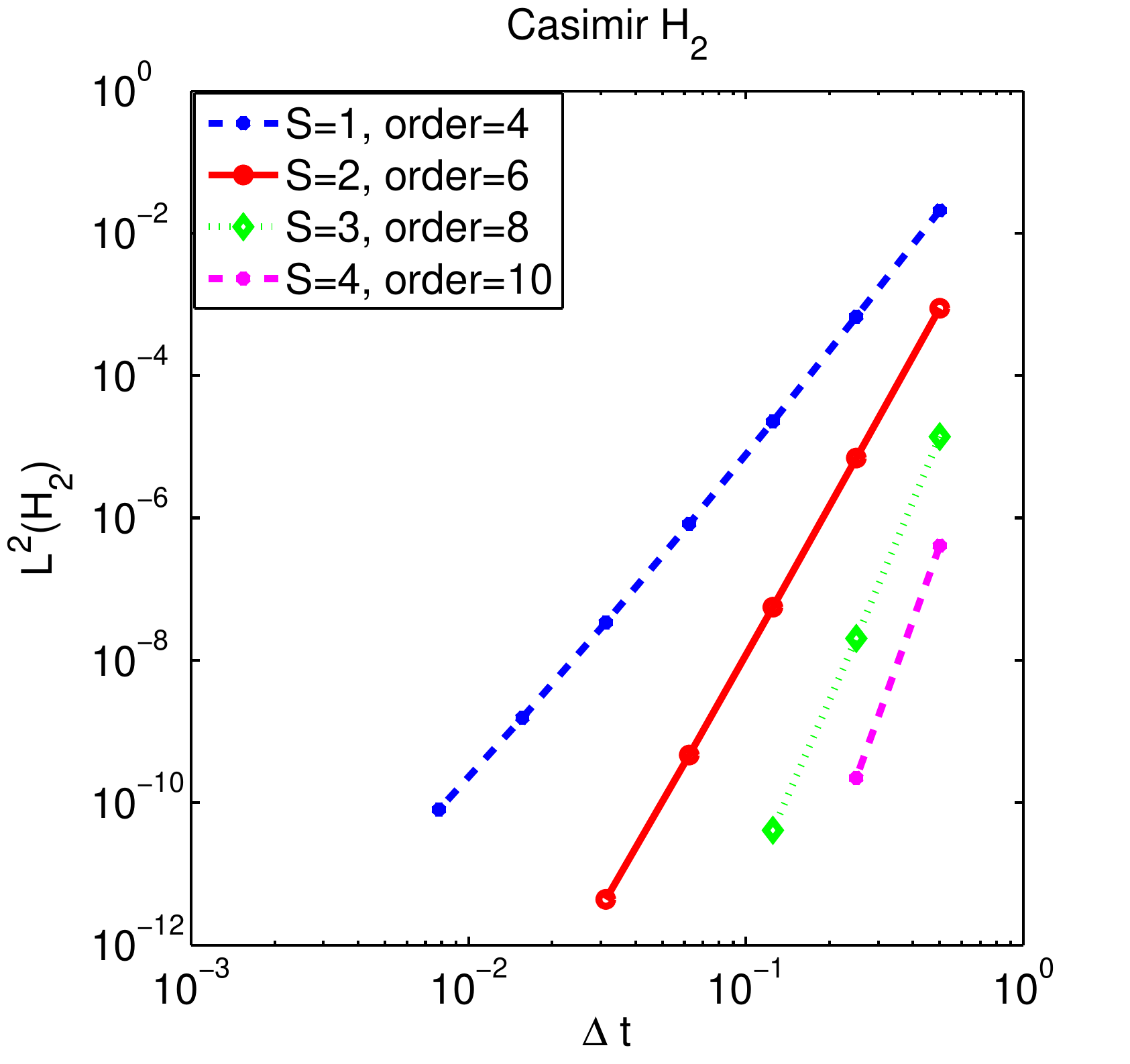}
\caption{Example~\ref{Ex1} with CDC: $L^2$-errors and convergence orders for the solution (left), Hamiltonians $H_1$ (right) and $H_2$ (right), with the choice $n=2S+3$}
\label{lv_fig_error}
\end{figure}

In Figure~\ref{lv_fig_error}, we give the $L^2$-errors and convergence orders of the solutions, Hamiltonians $H_1$ and $H_2$ for different correction number $S$  related with the number of nodes $n=2S+3$. When the errors reach about $10^{-10}$ level, the computations are stopped. With increasing number of correction step $S$, larger time steps are used to attain a prescribed order, which demonstrates the computational efficiency of CDC methods.

For varying correction number $S$ and node number $n$, the convergence orders are presented in Table~\ref{lv_table_mid} for the preservation of the Hamiltonian $H_1$. The results for the Hamiltonian $H_2$ and solutions are similar. The red labeled orders in Table~\ref{lv_table_mid} correspond to the setting $n=2S+3$, and they agree with the expected orders of accuracy. The time-step size $\Delta t$ in the table are chosen in order to reach the target error level $10^{-10}$ related to the runs with the red labeled orders.  The efficiency with respect to the step-size can be seen by the speedup factors in the last column, which are calculated as the ratio of the Wall Clock time required for the run without CDC method over the Wall Clock time required for the one with the CDC method, on the same level of accuracy $10^{-10}$.

\begin{table}[htb!]
\caption{Example~\ref{Ex1}: Convergence rates and speed-ups over the system without CDC \label{lv_table_mid}}
\centering
\begin{tabular}{c| c c c c c c c| c | c}
$S$/$n$ & 5 & 6 & 7 & 8 & 9 & 10 & 11 & $\Delta t$ & Speed-Up\\
\hline
 1  & \textcolor{red}{4.48}  & 4.53  & 4.65  & 4.47 & 4.75 & 4.72 & 4.83 & 0.01 & 3.0 \\
 2  & 4.51  & 6.99  & \textcolor{red}{6.78}  & 7.07 & 6.82 & 6.85 & 6.79 & 0.05 & 7.1 \\
 3  & 4.84  & 6.95  & 6.76  & 8.89 & \textcolor{red}{8.53} & 10.73 & 10.39 & 0.15 & 10.7 \\
 4  & 4.84  & 7.02  & 6.76  & 8.89 & 9.29 & 11.02 & \textcolor{red}{11.03} & 0.25 & 11.4 \\
\hline
\end{tabular}
\end{table}

%%%%%%%%%%%%%%%%%%%%%%%%%%%%%%%%%%%%%%%%%%%%%%%%%%%%%%%%%%%%%%%%%%%%%%%%%%%%%%%%%%%%%%%%
\subsection{Reversible LVS}
\label{Ex2}

\begin{figure}[htb!]
\centering
\includegraphics[scale=0.4]{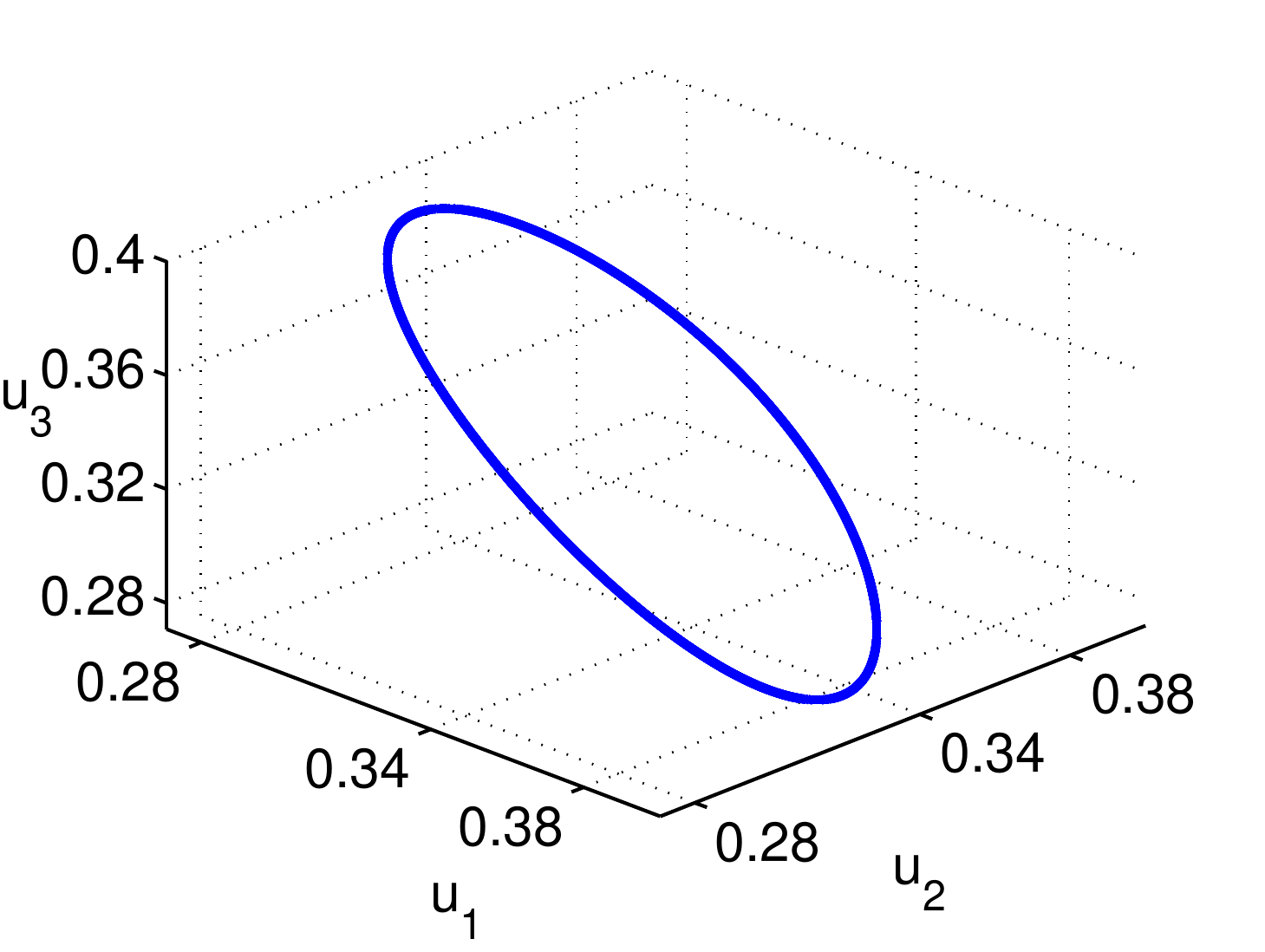}

\includegraphics[scale=0.4]{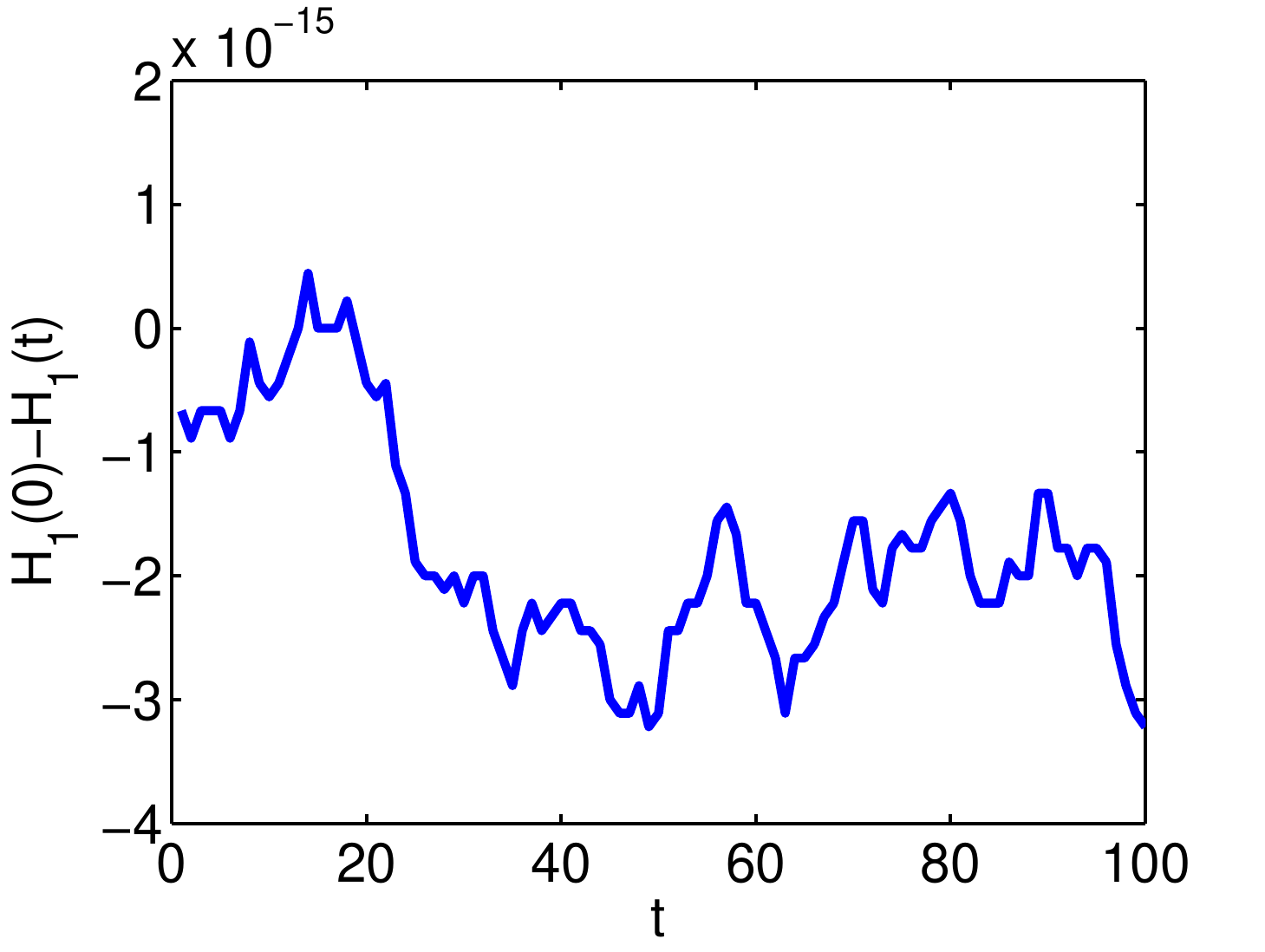}
\includegraphics[scale=0.4]{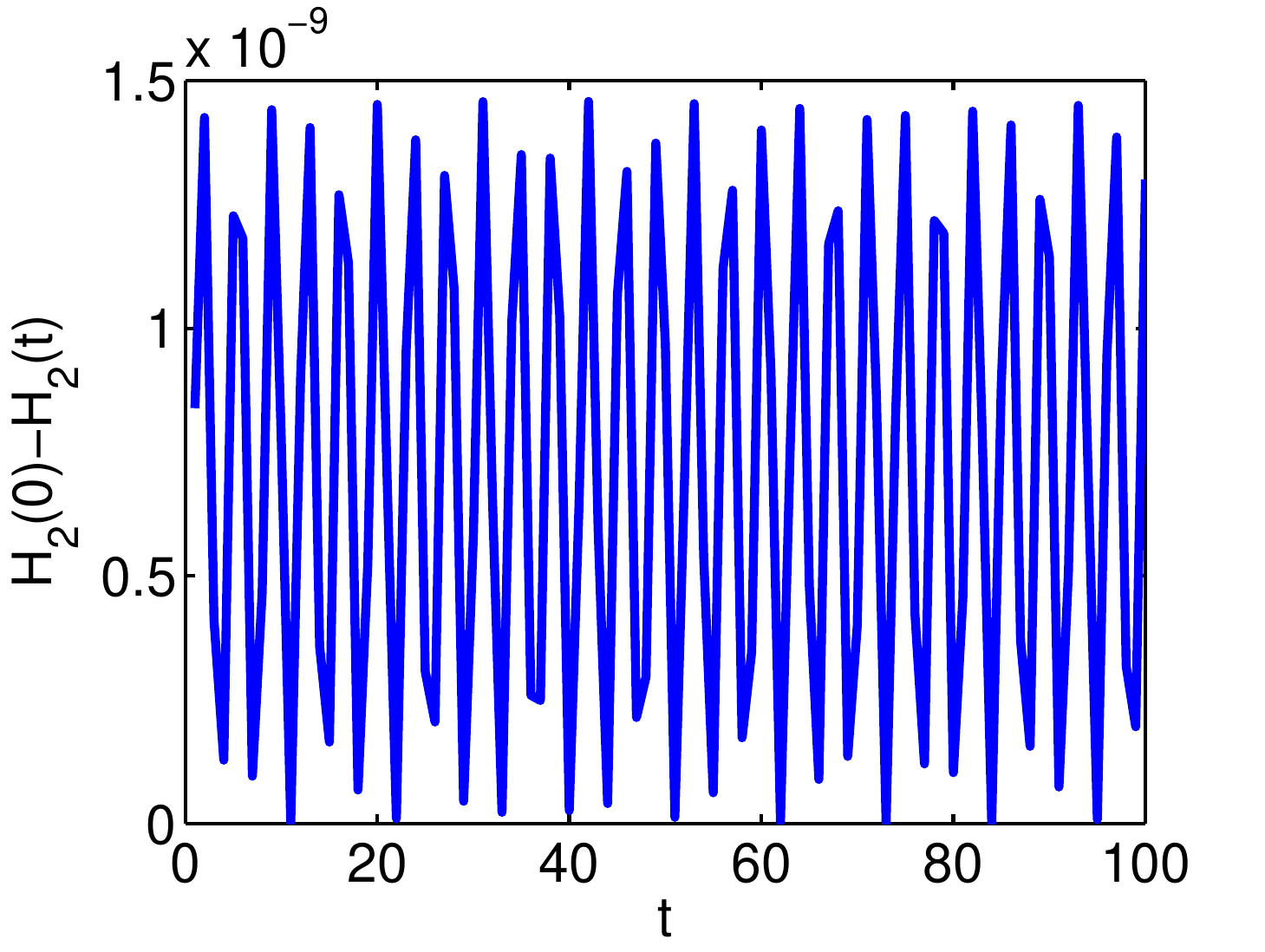}
\caption{Example~\ref{Ex2} without CDC: (Top) Periodic solutions; (Bottom) errors of Hamiltonians $H_1$ (left) and $H_2$ (right) from the initials: $\Delta t=0.01$}
\label{ps_fig_kahan}
\end{figure}

We consider the reversible LVS \eqref{lvs2} on the interval $[0,100]$, and with the initial conditions $(u_1(0),u_2(0),u_3(0) )^T = ( 0.3 , 0.3 , 0.4 )^T$ \cite{Ionescu15}. Kahan's method  preserves again the periodicity of the reversible LVS \eqref{lvs2}, Figure~\ref{ps_fig_kahan}, top.
The reversible LVS \eqref{lvs2} was solved in \cite{Wan17} using a conservative multiplier method. It was shown that the linear Hamiltonian $H_1$ is preserved with an accuracy $10^{-15}$ and the cubic Hamiltonian $H_2$ with an accuracy $10^{-14}$. Kahan's method also preserves the linear Hamiltonian $H_1$ and cubic Hamiltonian $H_2$  accurately in Figure~\ref{ps_fig_kahan}, bottom.

\begin{figure}[htb!]
\centering
\includegraphics[scale=0.34]{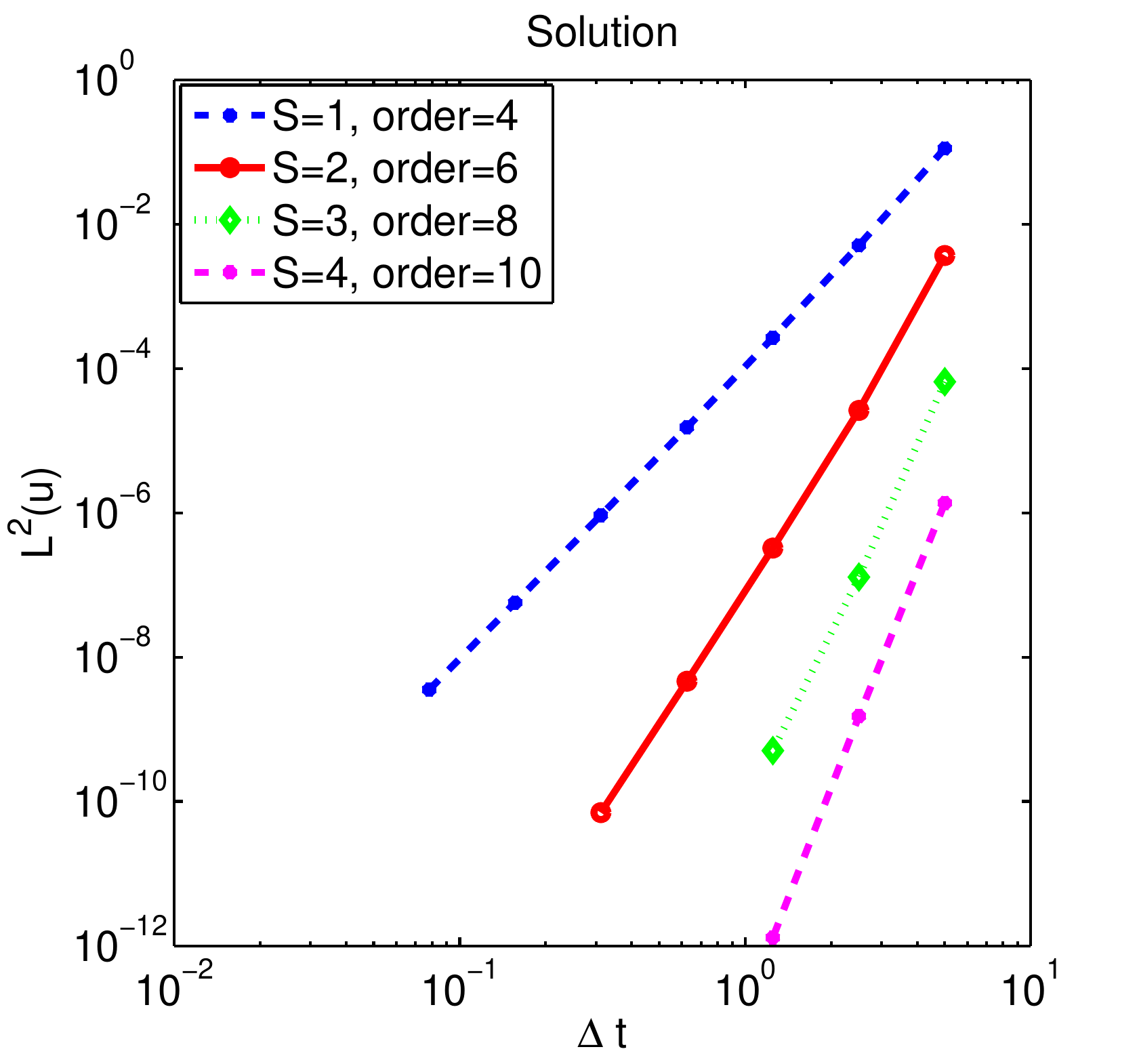}
\includegraphics[scale=0.34]{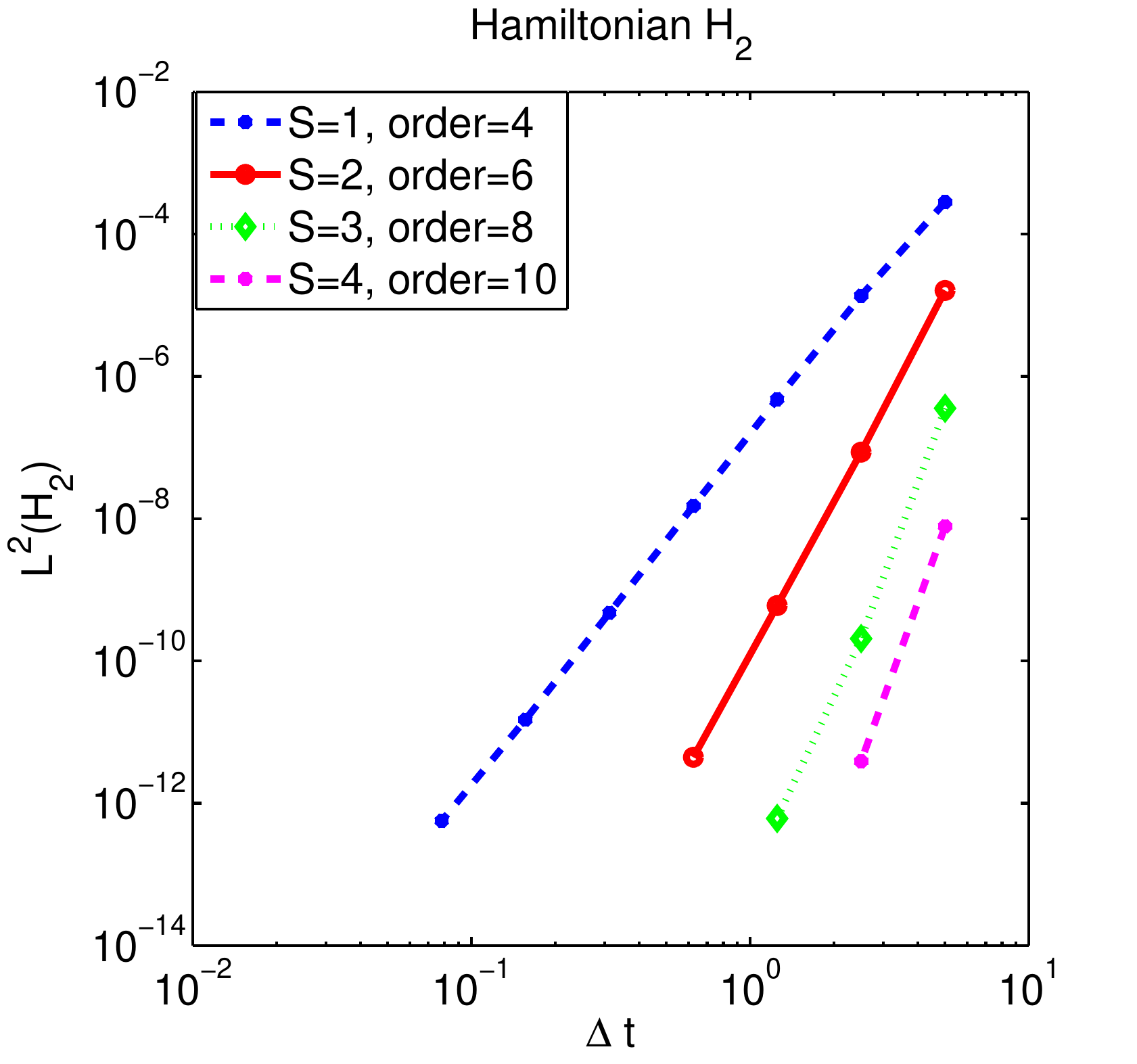}
\caption{Example~\ref{Ex2} with CDC: $L^2$-errors and convergence orders for the solution (left) and Hamiltonian $H_2$ (right), with the choice $n=2S+3$}
\label{ps_fig_error}
\end{figure}

Preservation of the periodicity of the solutions and the Hamiltonians in the case of CDC  method is similar to the LVS \eqref{lv1} in the previous example.
Again, similar convergence orders are attained in Figure~\ref{ps_fig_error}.

%%%%%%%%%%%%%%%%%%%%%%%%%%%%%%%%%%%%%%%%%%%%%%%%%%%%%%%%%%%%%%%%%%%%%%%%%%%%%%%
%%%%%%%%%%%%%%%%%%%%%%%%%%%%%%%%%%%%%%%%%%%%%%%%%%%%%%%%%%%%%%%%%%%%%%%%%%%%%%%
\section{Conclusions}
\label{sec:conclusion}

We have shown that the Hamiltonians of 3D LVSs can be preserved with a high accuracy, when we use CDC methods based on the Kahan's discretization for quadratic vector fields. In a future work, the integral and spectral correction methods on non-uniform grids will be applied. 

%%%%%%%%%%%%%%%%%%%%%%%%%%%%%%%%%%%%%%%%%%%%%%%%%%%%%%%%%%%%%%%%%%%%%%%%%%%%%%%

\end{document}